\documentclass[11pt]{article}
\usepackage{amsmath, amsthm, graphics, amssymb,fullpage, color,setspace, graphicx}

\newtheorem* {Main} {Theorem 1}

\newtheorem* {finitecenters} {Lemma 3.5}
\newtheorem* {prop1} {Proposition 3.1}
\newtheorem* {prop2} {Proposition 3.2}
\newtheorem* {permutecenters} {Lemma 3.6}
\newtheorem* {algind} {Lemma 3.3}
\newtheorem* {finitegroup} {Lemma 3.4}

\begin{document}

\title{Translation Surfaces With Finite Veech Groups}
\author{Asaf Hadari}
\maketitle
\begin{abstract} We show that every finite subgroup of $\textrm{GL}_{2}(\mathbb{R})$ can be realized as the Veech group of some translation surface.
\end{abstract}

\section{Introduction}

Veech groups play a pivotal role in the study of translation surfaces, quadratic differentials, and geodesics in the moduli space of Riemann surfaces. It is well known that Veech groups are always discrete subgroups of $\textrm{GL}_{2}(\mathbb{R})$, though it is not known which discrete subgroups of $\textrm{GL}_{2}(\mathbb{R})$ are Veech groups. A generic translation surface has a Veech group which is trivial or cyclic of order $2$. On the other end of the spectrum, some translation surfaces have Veech groups which are lattices. These are known as \textit{Veech surfaces}, and have been the object of much study.

The goal of this paper is to explore the smallest non-trivial possibilities: translation surfaces with finite Veech groups. The only finite subgroups of $GL_{2}(\mathbb{R})$ are cyclic and dihedral groups. We prove the following theorem.

\begin{Main}
Every finite subgroup of $GL_{2}(\mathbb{R})$ is can be realized as the  Veech group of some translation surface.
\end{Main}

Note that we allow both orientation-preserving and orientation-reversing elements in the Veech groups in this paper,that is, the Veech groups we consider are subgroups of $\textrm{GL}_{2}(\mathbb{R})$, and not just $\textrm{SL}_{2}(\mathbb{R})$.

Our method is constructive: we provide a translation surface for each such group. We make no effort at efficiency in terms of genus - for some finite groups there are examples of translation surfaces of lower genus that have the required Veech group.

\paragraph{Acknowledgements} The author wishes to thank Benson Farb ,Howard Masur, Alex Eskin, and Matthew Bainbridge for useful discussions on the ideas in the paper. He also wishes to thank his wife Nurit Kirshenbaum for creating  the figures.

\section{Preliminaries}
\paragraph{Translation Surfaces.}
A translation surface $T$ is a $2$-dimensional manifold containing a discrete subset $\Sigma \subset T$ such that $T \backslash \Sigma$ is equipped with a maximal atlas with the property that the transition functions are translations. The set $\Sigma$ is called the set of \textit{cone points} of $T$. Note that the atlas above imbues $T$ with a flat metric away from the set of cone points.

One way to construct translation surfaces is the following: start with a polygon in $\mathbb{R}^{2}$ with the property that each side of the polygon is parallel and congruent to a different edge. By gluing such edges in pairs, one obtains a translation surface. For example, by gluing parallel congruent edges of a rectangle in the plane one obtains a flat torus with no cone points. By gluing parallel congruent edges of a $4g$-gon ($g \geq 2$) one obtains a genus $g$ surface with one cone point.

Note that there are several equivalent ways in the literature of defining translation surfaces. We choose the point of view which is simplest for our needs.

\paragraph{The Developing Map and Holonomy.} Let $T$ be a translation surface with cone points $\Sigma$. Let $\tilde{T}$ be the universal cover of $T \backslash \Sigma$. The manifold $\tilde{T}$ has a flat metric, given by pulling back the metric from $T \backslash \Sigma$. Given any choice of basepoint $p \in \tilde{T}$, there is a unique locally isometric embedding $\tilde{T} \to \mathbb{R}^{2}$ sending $p$ to the origin in $\mathbb{R}^{2}$. This map is called the \textit{developing map}. Given a path $\gamma \subset T \backslash \Sigma$ and a lift $\tilde{\gamma}$ of $\gamma$ to $\tilde{T}$, the difference of the image of the endpoints of $\tilde{\gamma}$ is independent of the lift and of the choice of basepoint $p$. We denote this difference $\textrm{hol}(\gamma)$. If $T$ itself is simply connected and a basepoint $p$ is chosen, then the holonomy map is a well defined map from $T$ to $\mathbb{R}^{2}$. We denote this map  $\textrm{hol}_{p}$.

\paragraph{Saddle Connections.} A \textit{saddle connection} is a geodesic segment whose endpoints are cone points, and whose interior does not contain any cone points. Notice that since the set of cone points is discrete, there are only finitely many saddle connections on a given translation surface whose lengths are less than or equal to a given number.

One important fact that we use about saddle connections is the following: let $\textrm{Hol}(T)$ be the $\mathbb{Z}$-module generated by the holonomy of all saddle connections. Then $\textrm{Hol}(T)$ is finitely generated as a $\mathbb{Z}$-module (\cite{KeSm}).

\paragraph{Affine Groups and Veech Groups.} Given a translation surface $T$ with set of cone points $\Sigma$, the flat structure gives a trivialization of the tangent bundle of $T \backslash \Sigma$. Thus, to each diffeomorphism $f: T \backslash \Sigma \to T \backslash \Sigma$ one can attach a map $T \backslash \Sigma \to \textrm{GL}_{2}(\mathbb{R})$ which assigns to a point $p$ the derivative of $f$ at $p$.  The group consisting of all maps $f: T \to T$ that fix $\Sigma$, and that have constant derivatives away from $\Sigma$ is called the \textit{affine group of T}, and we denote it $\textrm{Aff}(T)$. The group $\textrm{Aff}(T)$ has a natural projection $\textrm{Aff}(T) \to \textrm{GL}_{2}(\mathbb{R})$ given by taking derivatives. We denote its image by $\textrm{GL}_{2}(T)$. This group is called the \textit{Veech group of T}.

Note that the group $\textrm{GL}_{2}(T)$ is actually a subgroup of $\textrm{SL}_{2}^{\pm}(\mathbb{R})$, where $\textrm{SL}_{2}^{\pm}(\mathbb{R})$ is the group of all elements whose determinant is $\pm 1$. To see this, consider the module $\textrm{Hol}(T)$ defined above. Since $\textrm{Aff}(T)(\Sigma) = \Sigma$, and elements of $\textrm{Aff}(T)$ send segments to segments, one has that $\textrm{GL}_{2}(T) \textrm{Hol}(T) = \textrm{Hol}(T)$. It's a standard fact that the stabilizer of any finitely generated $\mathbb{Z}$-submodule of $\mathbb{R}^{2}$ must be a subgroup of $\textrm{SL}_{2}^{\pm}(\mathbb{R})$.

Recall that an element of $SL_{2}^{\pm}(\mathbb{R})$ is \textit{elliptic} if its trace has absolute value $< 2$, \textit{parabolic} if its trace has absolute value $2$, and \textit{hyperbolic} if its trace has absolute value $> 2$. We can thus classify all elements of the Veech group of a translation surface as elliptic, parabolic, or hyperbolic. In this paper we are predominantly interested in elliptic elements of Veech groups. These always have finite order, and preserve the flat metric on the translation surface.

\section{Proof of Theorem $1$.}

 Any finite subgroup $G < GL_{2}(\mathbb{R})$ preserves an Euclidean metric on $\mathbb{R}^{2}$, and can thus be conjugated into $O_{2}(\mathbb{R})$. The finite subgroups of $O_{2}(\mathbb{R})$ are exactly the finite cyclic groups and the finite dihedral groups. We denote the cyclic group of order $N$ by $C_{N}$. We denote the dihedral group of order $2N$ by $D_{N}$. In order to prove Theorem $1$, we prove the following two propositions.

\begin{prop1} For every $N \geq 3$, there exists a translation surface $T = T(N)$ such that $\textrm{GL}_{2}(T) \cong D_{N}$.
\end{prop1}

\begin{prop2} For every $N \geq 3$, there exists a translation surface $T = T(N)$ such that $\textrm{GL}_{2}(T) \cong C_{N}$.
\end{prop2}

The proofs of the the two propositions are very similar. The following lemmas are used in the proof of both propositions.

\begin{algind} Suppose $T$ is a translation surface such that $\textrm{Hol}(T)$ contains two $\mathbb{R}$-linearly independent and algebraically independent vectors. Then $\textrm{GL}_{2}(T)$ contains no hyperbolic elements.

\end{algind}
\paragraph{Proof.}
Let $e_{1}$ and $e_{2}$ be two vectors as in the conditions of the lemma. Suppose that there was a hyperbolic element $T \in \textrm{GL}_{2}(T)$ with positive trace. Note that the square of any hyperbolic element is always a hyperbolic element with positive trace.  Let $K = \mathbb{Q}[\textrm{Trace}(T)]$. The field $K$ is a number field, and $\textrm{Hol}(T) \otimes \mathbb{Q}$ is a
 $2$-dimensional vector space over $K$ (see \cite{KeSm} or \cite{MC1} for a different proof). However, since $e_{1}, e_{2} \in \textrm{Hol}(T)$ are algebraically independent, we must have that $\textrm{Hol}(T) \otimes \mathbb{Q}$ is has infinite dimension as a $\mathbb{Q}$, and thus as a $K$-vector space. This is a contradiction, thus $\textrm{PGL}_{2}(T)$ contains no hyperbolic elements. $\Box$

\begin{finitegroup}
Let $T$ be a translation surface such that $\textrm{GL}_{2}(T)$ contains no hyperbolic elements. Suppose $\textrm{GL}_{2}(T)$ contains an elliptic element of order at least $3$. Then $\textrm{GL}_{2}(T)$ is finite.
\end{finitegroup}

\paragraph{Proof.} Suppose that $\textrm{GL}_{2}(T)$ contained a parabolic element with positive trace (once again, note that the square of any parabolic element is always a parabolic element of positive trace.) We denote this element by $P$. Let $E$ be an elliptic element of  $\textrm{GL}_{2}(T)$ of order at least $3$. Let $Q = EPE^{-1}$. The transformation $Q$ is a parabolic element of the Veech group that does not commute with $P$. By a construction due to Thurston, a Veech group that contains two non-commuting parabolic elements also contains hyperbolic elements (See \cite{Thurston} for the construction, or \cite{Veechcalc} for an application of this construction to Veech groups.) This contradicts our assumption.

Therefore, $GL_{2}(T)$ is composed entirely of elliptic elements. Any discrete subgroup of $\textrm{SL}_{2}(\mathbb{R})$ whose elements are all elliptic is finite (\cite{katok}, page $37$). Any subgroup of $\textrm{SL}_{2}^{\pm}(\mathbb{R})$ contains a finite (at most $2$) index subgroup contained in $\textrm{SL}_{2}(\mathbb{R})$. Thus $\textrm{GL}_{2}(T)$ is finite.  $\Box$

\bigskip

Another tool that is used several times in the proof of Theorem $1$ is a gadget we call "maximum embedded convex polygons", which exploits the fact that elliptic elements are isometries in the flat metric in order to show that they must permute certain finite sets.

\paragraph{Definition.} Let $T$ be a translation surface. Let $\Sigma$ be its set of cone points. A proper closed subset $C \subset T$ is a \textit{maximal embedded convex polygon in T} if:

\begin{enumerate}
\item $\pi_{1}(C^{o}) = 1$, where $C^{o}$ denotes the interior of $C$.
\item $\exists p \in C^{0}$ such that $\textrm{hol}_{p}(C)$ is a closed convex polygonal region in $\mathbb{R}^{2}$.
\item $C^{o} \cap \Sigma = \emptyset$
\item $\partial C$ is a finite union of saddle connections.
\end{enumerate}

\paragraph{Definition.} Let $T$ be a translation surface. Let $C \subset T$ be a maximal embedded convex polygon, and $p \in C$. Let $\hat{c}$ be the Euclidean center of mass of $\textrm{hol}_{p}(C)$. Since the developing map is injective on $C^{o}$, there is a unique point $c: = \textrm{hol}_{p}^{-1}(\hat{c})$. We call $c$ the \textit{centroid of C}. Note that this point is independent of the choice of $p$.

\paragraph{Definition.} Let $T$ be a translation surface. Let $P \subset \mathbb{R}^{2}$ be a closed convex polygonal region. We define $\textrm{Copies}(T:P)$ to be the set of all maximal embedded convex polygons $C \subset T$ such that $\textrm{hol}_{p}(C)$ is isometric to $P$ for any $p \in C$.

\paragraph{Definition.} Let $T$ be a translation surface. Let $P \subset \mathbb{R}^{2}$ be a closed convex polygonal region. Let $\textrm{Centroids}(T:P)$ equal the set of centroids of all the elements of $\textrm{Copies}(T:P)$.

\begin{finitecenters}
Given a translation surface $T$, and a closed convex polygonal region $P \subset \mathbb{R}^{2}$, the sets $\textrm{Copies}(T:P)$ and $\textrm{Centroids}(T:P)$ are finite.
\end{finitecenters}

\paragraph{Proof.} Let $M = \textrm{diameter}(P)$. All of the sides of $\partial P$ have length $\leq M$. Let $C \in \textrm{Copies}(T:P)$. The set $\partial C$ is composed of saddle connections of length $\leq M$. The interior $C^{o}$ is  a connected component of $T \backslash \partial C$. Since there are only finitely many saddle connections of length $\leq M$, and only finitely many connected components obtained by removing a subset of them from $T$, we conclude that $\textrm{Copies}(T:P)$ is finite. Since $\# \textrm{Copies}(T:P) \geq \# \textrm{Centroids}(T:P)$, we conclude that $\textrm{Centroids}(T:P)$ is finite as well.  $\Box$

\begin{permutecenters}
Given a translation surface $T$, a closed convex polygonal region $P \subset \mathbb{R}^{2}$, and an elliptic element $E \in \textrm{Aff}(T)$, the following holds:

$$E(\textrm{Centroids}(T:P)) = \textrm{Centroids}(T:P) $$
\end{permutecenters}

\paragraph{Proof.} The map $E$ is an isometry of the flat metric. Furthermore, it preserves the set of cone points, and thus sends saddle connections to saddle connections. Thus,  $E(\textrm{Copies}(T:P) = \textrm{Copies}(T:P)$.
Given $C_{1}, C_{2} \in \textrm{Copies}(T:P)$ with $E(C_{1}) = C_{2}$ and $p \in C_{1}$,  then $\textrm{hol}_{E(p)} \circ E \circ \textrm{hol}_{p}^{-1}$ is an isometry from $\textrm{hol}_{p_{1}}(C_{1})$ to $\textrm{hol}_{E(p_{1})}(C_{2})$. Thus, $E$ must send the centroid of $C_{1}$ to the centroid of $C_{2}$. $\Box$

\paragraph{Proof of Proposition $3.1$.}

We split into two cases - one where $N$ is odd, and the other where $N$ is even. Most of the important features of the proofs of Proposition $3.1$, and Proposition $3.2$ appear in the case where $N$ is even, so we go over this case in the most detail.

\paragraph{The case where $N$ is even.} Choose numbers $0 < \ell_{2} < \frac{1}{3} \ell_{1} < \ell_{1}$, such that $\ell_{1}$ and $\ell_{2}$ are algebraically independent. Construct an isosceles  triangle whose congruent sides meet at the point $O$ at an angle of $\theta: = \frac{2 \pi}{N}$, such that the side opposite $O$ has length $\ell_{1}$. Construct a square with sides of length $\ell_{2}$ such that one side of the square is contained in the side of the triangle opposite $O$, and such that the midpoints of these two sides coincide. Call the center of this square $P_{1}$. Rotate the resulting figure $N - 1$ times about the point $O$ by angles $\theta, 2 \theta, \ldots, (N -1) \theta$. Label the rotation images of $P_{1}$: $P_{2}, \ldots P_{N}$. The resulting figure is a polygon in the plane, consisting of a regular $N$-gon, with $N$ squares attached to the midpoints of its sides.

 Choose a pairing of each side of the polygon with a parallel congruent side such that the pairing is preserved by rotating the polygon about $O$ by an angle of $\theta$, and by reflecting over a line of symmetry of the polygon. Identifying the edges paired in this way, we obtain  a translation surface $T$. The cone points of $T$ are precisely the images under this identification of all of the vertices of the polygon which lie on the regular $N$-gon. The construction is illustrated in Figure $1$ for the case $N = 4$.

\begin{figure} [h!]
\begin{center}
\includegraphics[width=0.5\textwidth]{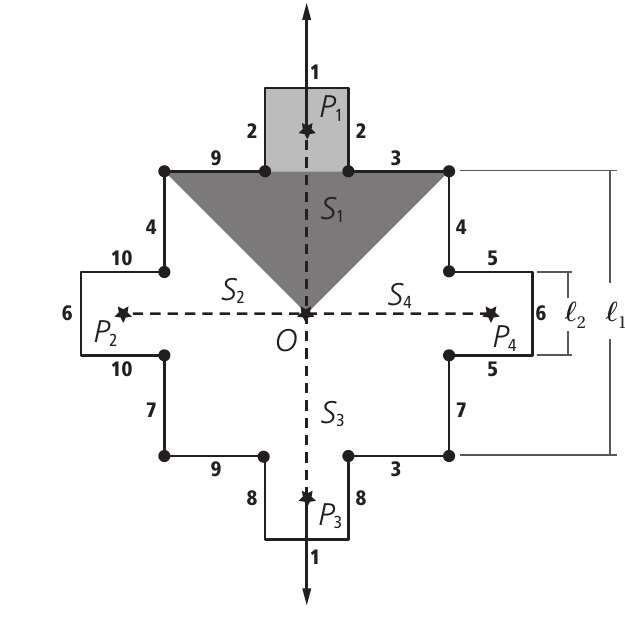}
\caption{The planar polygon and the edge identification for the case $D_{4}$.}
\end{center}
\end{figure}

Notice that $\textrm{Aff}(T)$ contains a subgroup isomorphic to $D_{N}$. Indeed, the polygon we constructed in the plane clearly has such a subgroup of symmetries (given by rotating it about $O$ by integer multiples of $\theta$ and by reflecting over the line $\mathcal{L}$ pictured in the figure.) Since these symmetries respect the identifications of parallel edges, we get affine maps of $T$. Since none of the above symmetries has trivial derivative, we get that $\textrm{GL}_{2}(T)$ contains a subgroup isomorphic to $D_{N}$.

\bigskip
Our next goal is to show that $\textrm{GL}_{2}(T) \cong D_{N}$. In order to do this, we must show that there are no other elements in $\textrm{Aff}(T)$ aside from the ones described above.

We first notice that $\textrm{GL}_{2}(T)$ is finite. Indeed, $\ell_{1} - \ell_{2}, i\ell_{2} \in \textrm{Hol}(T)$, and thus by Lemma $3.3$, $\textrm{GL}_{2}(T)$ contains no hyperbolic elements. Since it contains an elliptic element of order $\geq 3$, we have by Lemma $3.4$ that $\textrm{GL}_{2}(T)$ is finite.

Let $Q_{1} \subset \mathbb{R}^{2}$ be a regular $N$-gon with sides of length $\ell_{1}$. We have that $\# \textrm{Copies}(T:Q_{1}) = 1$, and $\textrm{Centroids}(T:Q_{1}) = O$. One can see this by checking that $O$ is the only point in $T$ whose distance from the set of cone points is $\frac{\ell_{2}}{\sin (\frac{\theta}{2})}$. All other points are closer to the set of cone points. Thus, by Lemma $3.6$, the point $O$ is $\textrm{Aff}(T)$-invariant. Notice that we now have that $\textrm{Aff}(T) \cong \textrm{Gl}_{2}(T)$. Indeed, suppose $A \in \textrm{Aff}(T)$, with $dA = Id$. Since $A$ is an isometry, $A(O) = O$, and $dA = Id$, we have that $A = Id$. In what follows we will make no distinction between elements of $\textrm{Aff}(T)$ and their images in $\textrm{GL}_{2}(T)$.

Let $Q_{2} \subset \mathbb{R}^{2}$ be a square with sides of length $\ell_{2}$. It is simple to see that $$\textrm{Centroids}(T:Q_{2}) = \{P_{1}, \ldots P_{N} \}.$$ Indeed, $\ell_{2}$ is the length of the shortest saddle connection in $T$, and the collection of these saddle connections forms $N$ squares with centroids $\{P_{1}, \ldots P_{N} \}$. Thus, by Lemma $3.6$, the set $\{P_{1}, \ldots P_{N} \}$ is $\textrm{GL}_{2}(T)$ invariant.

Given any $i = 1, \ldots N$, there is a unique geodesic segment in $T$ of length $\frac{\ell_{1}}{2 \tan{\frac{\theta}{2}}} + \frac{\ell_{2}}{2}$ connecting $O$ to $P_{i}$, and containing no cone points. Indeed, each $P_{i}$ is the centroid of an element of $\textrm{Copies}(T:Q_{2})$, and $O$ is the centroid of the unique element in $\textrm{Copies}(T:Q_{1})$. The distance of $O$ from the boundary of the regular $N$-gon isometric to $Q_{1}$ is $\frac{\ell_{1}}{2 \tan{\frac{\theta}{2}}}$, and this distance is realized exactly by the line segments connecting $O$ to the midpoints of the sides of the $N$-gon. Similarly, the distance of $P_{1}$ from the boundary of the square isometric to $Q_{2}$ is   $\frac{\ell_{2}}{2}$, and this distance is realized exactly by the segments connecting $P_{i}$ to the midpoints of one of the sides of the square. Since the interior of the regular $N$-gon is disjoint from the interiors of all of the squares, we get that each such segment must connect $O$ to the midpoint of one of the sides of the $N$-gon, and then connect that midpoint to $P_{i}$. There are exactly $N$ such segments, which we label $S_{1}, \ldots, S_{N}$. The segment $S_{2}$ is shown in Figure $1$.
\bigskip

Since the sets $\{O \}$ and $\{P_{1}, \ldots P_{N} \}$ are $\textrm{GL}_{2}(T)$ invariant, then the set of segments joining $O$ to $\{P_{1}, \ldots P_{N} \}$ is $\textrm{GL}_{2}(T)$ invariant. Since $\textrm{GL}_{2}(T)$  acts on $T$ by isometries, we get that the set $\{S_{1}, \ldots S_{N} \}$ is $\textrm{GL}_{2}(T)$ invariant. Let $\sigma_{1}, \ldots \sigma_{N} \in T_{O}(T)$ be elements of the tangent space at $O$ such that $\sigma_{i} = \textrm{hol}(S_{i})$. That is - $\sigma_{i}$ is a vector in the tangent space to $T$ at $O$ with the same length and direction as $S_{i}$. The group $\textrm{GL}_{2}(T)$ fixes $O$, and thus acts on $T_{O}(T)$ by isometries. The set $\{\sigma_{1}, \ldots, \sigma_{N}\}$ is clearly $\textrm{GL}_{2}(T)$-invariant.

Consider the derivative representation $\pi: \textrm{GL}_{2}(T) \to \textrm{GL}_{2}(T_{O}(T))$. Notice that $\textrm{Ker}(\pi) = Id$.
We claim that $\pi(\textrm{GL}_{2}(T)) \cong D_{N}$. Indeed, we already showed that $\pi(\textrm{GL}_{2}(T))$ contains a subgroup isomorphic to $D_{N}$. Since  every element of $\textrm{Gl}_{2}(T)$ is an isometry, we must have that that $\pi(\textrm{GL}_{2}(T)) \subset O_{2}(\mathbb{R})$. Furthermore, the set $\{\sigma_{1}, \ldots \sigma_{N} \}$ is the set of vertices of a regular $N$-gon with centroid at the origin. Since this set is $\pi(\textrm{GL}_{2}(T))$-invariant, and the $O_{2}(\mathbb{R})$-stabilizer of such a set is isomorphic to $D_{N}$, we must have $\pi(\textrm{GL}_{2}(T)) \cong D_{N}$. Since $\pi$ is injective, this means that $\textrm{GL}_{2}(T) \cong D_{N}$.

\paragraph{The case where $N$ is odd.} We first note that the construction detailed in the previous case does not work for odd values of $N$. Indeed, the polygon constructed in the previous case has no parallel sides when $N$ is odd, and thus cannot be used to construct a translation surface as before. We rely on a slightly more complicated construction instead. The reader may note that this more complicated construction works for even values of $N$ as well.

\bigskip

We begin similarly to the previous construction. Choose numbers $0 < \ell_{2} < \frac{1}{3} \ell_{1} < \ell_{1}$, such that $\ell_{1}$ and $\ell_{2}$ are algebraically independent. Construct an isosceles  triangle whose congruent sides meet at the point $O$ at an angle of $\theta: = \frac{\pi}{N}$, such that the side opposite $O$ has length $\ell_{1}$. Construct a square with sides of length $\ell_{2}$ such that one side of the square is contained in the side of the triangle opposite $O$, and such that the midpoints of these two sides coincide. Call the center of this square $P_{1}$. The square has two sides which are perpendicular to the base of the isosceles triangle. The midpoints of these sides divide them into four segments: $A_{1}, \ldots, A_{4}$, as illustrated in the Figure $2$.

\begin{figure} [h!]
\begin{center}
\includegraphics[width=0.5\textwidth]{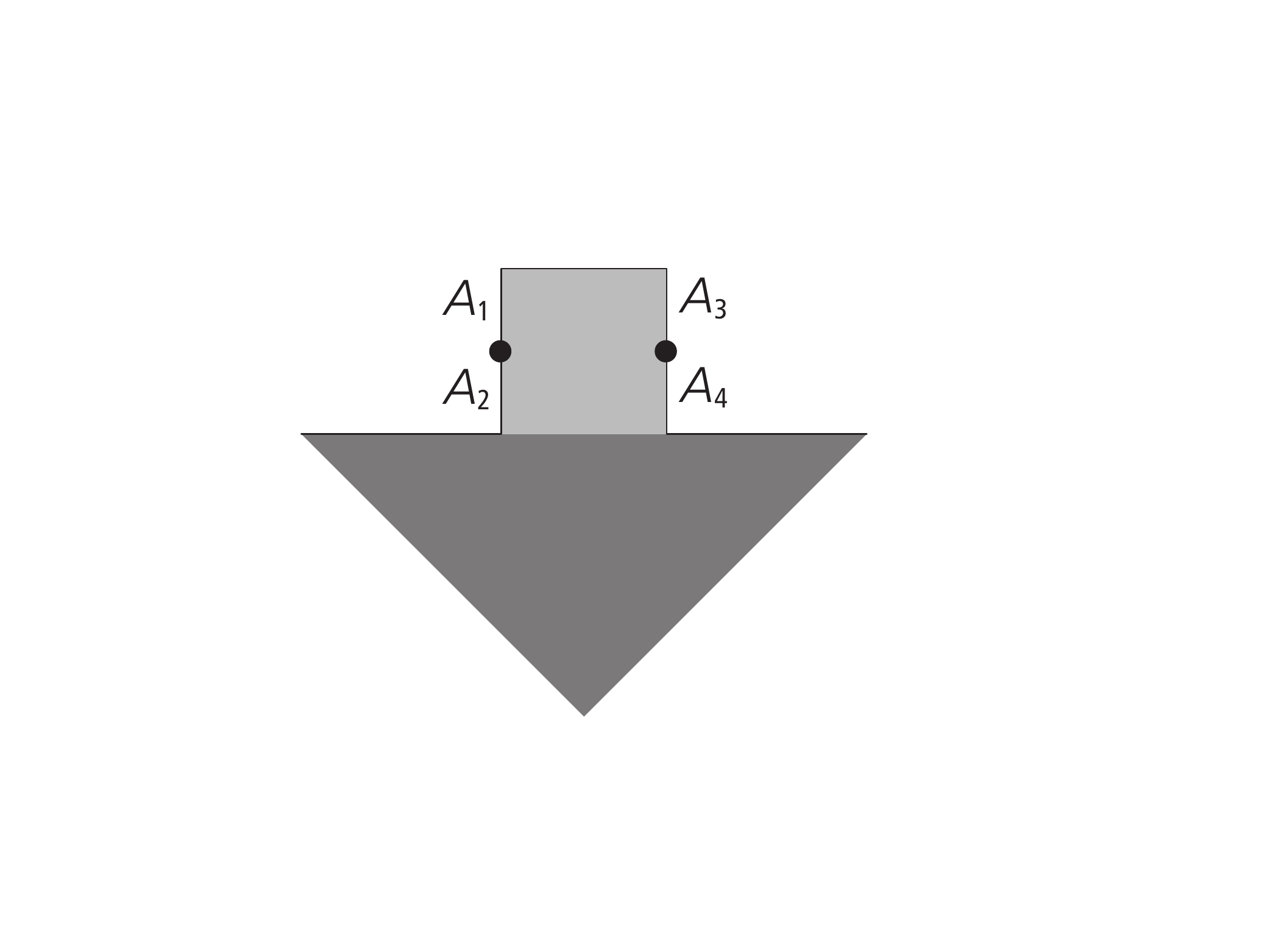}
\caption{The segments $A_{1}, A_{2}, A_{3}, A_{4} $.}
\end{center}
\end{figure}

Rotate the resulting figure $2N -1$ times about the point $O$ by angles $\theta, 2 \theta, \ldots, (2N - 1) \theta$. Label the rotation images of $P_{1}$: $P_{2}, \ldots P_{2N}$. For $i = 1, \ldots, 4$ and $j = 1, \ldots, 2N$, label the rotation image of $A_{i}$ by $A_{i + 4j}$.

The resulting figure is a polygon in the plane. We now wish to identify parallel edges to obtain a translation surface. For $j = 1, 3, \ldots, 2N - 1$ identify $A_{1 + 4j}$ with $A_{4 + 4j}$, and $A_{2+4j}$ with $A_{3+4j}$. For $j = 0, 2, 4, \ldots 2N -2$, identify $A_{1+4j}$ with $A_{3+4j}$ and $A_{2+4j}$ with $A_{4+4j}$. For all other sides of the polygon, identify each side with one of the sides which is parallel and congruent to it, in a manner which is invariant under rotation by an angle of $\theta$ and under reflection. Call the resulting translation surface $T$. This construction is illustrated for the case $N = 3$ in Figure $3$.

\begin{figure} [h!]
\begin{center}
\includegraphics[width=0.5\textwidth]{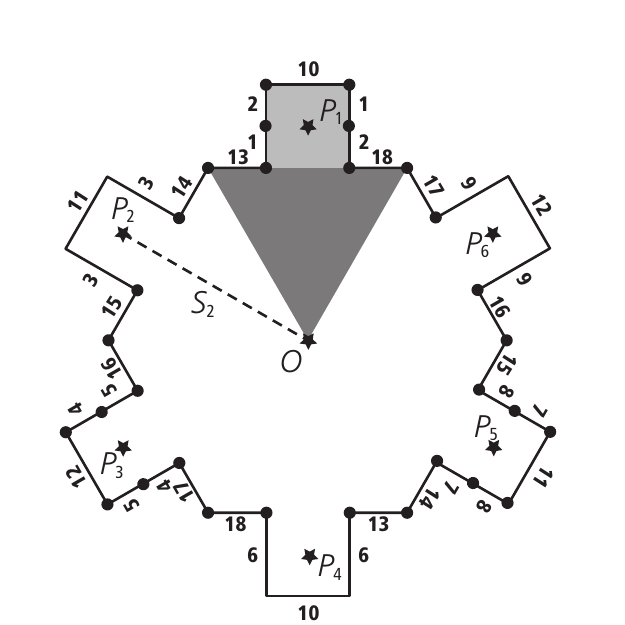}
\caption{The polygon and edge identifications for the case $D_{3}$.}
\end{center}
\end{figure}

The proof now proceeds similarly to the previous case. By the same considerations as the ones in the previous case, $\textrm{Aff}(T)$ is a finite group consisting of elliptic elements. Furthermore, this group is isomorphic to $\textrm{GL}_{2}(T)$, fixes the point $O$, permutes the set of points $\{P_{1}, \ldots, P_{2N}\}$, and the set of segments $\{S_{1}, \ldots, S_{2N}\}$ which were defined in the previous case.

Note that $\textrm{GL}_{2}(T)$ contains a subgroup isomorphic to $D_{N}$ which is generated by the reflection over a line of symmetry of the polygon, and a rotation about $O$ by an angle of $2\theta$. We now wish to show that there are no other elements in $\textrm{GL}_{2}(T)$.

Let $\Sigma$ be the set of cone points of $T$. This set consists of the images under the identification above of the following points: the vertices of the polygon that lie on the regular $N$-gon, and the endpoints of $A_{1 + 4j}, \ldots, A_{4+4j}$ for $j$ odd.

Suppose $i$ is an even number. It's clear that $\textrm{dist}(P_{i}, \Sigma) = \frac{\sqrt{2}}{2} \ell_{2}$. On the other hand, for $i$ odd, we have that $\textrm{dist}(P_{i}, \Sigma) = \frac{\ell_{2}}{2}$. Suppose we are given an odd index $i_{o}$ and an even index $i_{e}$, and  $E \in \textrm{GL}_{2}(T)$. Since $E$ is an isometry and $E(\Sigma) = \Sigma$, we cannot have that $E(P_{i_{o}}) = P_{i_{e}}$. Similarly, it is impossible to have $E(S_{i_{o}}) = S_{i_{e}}$. Thus, the set $\{S_{1}, S_{3}, S_{5}, \ldots, S_{2N - 1} \}$ is $\textrm{GL}_{2}(T)$-invariant. Proceeding as in the previous case, we see that $\textrm{GL}_{2}(T) < D_{N}$, and thus $\textrm{GL}_{2}(T) \cong D_{N}$. $\Box$

\paragraph{Proof of Proposition $3.2$.} Once again we divide into two cases, depending on the parity of $N$. The proofs are very similar to the two cases in the proof of Proposition $3.1$, and we omit details of the proof when they are nearly identical to those found in the proof of Proposition $3.1$.

\paragraph{The case where $N$ is even.} Choose numbers $0 < \ell_{3} < \frac{1}{2}\ell_{2} <  \ell_{2} < \frac{1}{3} \ell_{1} < \ell_{1}$, such that $\ell_{1}$, $\ell_{2}$ and $\ell_{3}$ are algebraically independent. Construct a polygon as in the even case of the proof of Proposition $3.1$, with the following difference: when constructing the square with sides of length $\ell_{2}$ that shares a side with the isosceles triangle, we require that the midpoints of the sides of the square and the triangle \textit{do not} coincide. Let the distance between the midpoints be $\ell_{3}$. As in the the previous proofs, we get a polygon in the plane. We label the points $P_{1}, \ldots, P_{N}$ as before. Furthermore, for each $1 \leq i \leq N$, there is a unique segment connecting $O$ to $P_{i}$. As in the previous cases, we label this segment $S_{i}$. Label the $N$ segments connecting $O$ to the vertices of the regular $N$-gon: $T_{1}, \ldots, T_{N}$.

Identifying each edge with the unique edge which is parallel and congruent to it, we get a translation surface $T$. The construction is illustrated in Figure $4$ for the case $N = 4$.

\begin{figure} [h!]
\begin{center}
\includegraphics[width=0.5\textwidth]{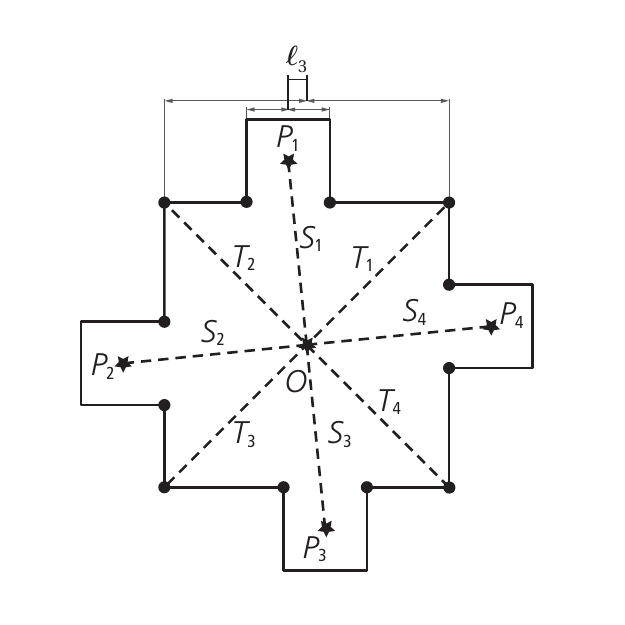}
\caption{The polygon for the case $C_{4}$.}
\end{center}
\end{figure}

\bigskip

By the same considerations as the ones in the previous cases, $\textrm{Aff}(T)$ is a finite group consisting of elliptic elements. Furthermore, this group is isomorphic to $\textrm{GL}_{2}(T)$, fixes the point $O$, permutes the set of points $\{P_{1}, \ldots, P_{N}\}$, and the set of segments $\{S_{1}, \ldots, S_{N}\}$.

Notice that each of the points in $T$ corresponding to a vertex of the regular $N$-gon is a cone point. Each of these points has distance $\frac{\ell_{1}}{2\sin\frac{\theta}{2}}$ from $O$, and that $T_{1}, \ldots, T_{N}$ are the only segments of this length connecting $O$ to the set of cone points. Thus, $\textrm{GL}_{2}(T)$ permutes the set $\{T_{1}, \ldots, T_{N}\}$.

Notice that $\textrm{GL}_{2}(T)$ contains a subgroup isomorphic to $C_{N}$. These elements come from rotating the polygon about $O$ by an angle of $\theta$. As before, we now wish to show that there are no other elements in $\textrm{GL}_{2}(T)$.

Consider the derivative representation $\pi: \textrm{GL}_{2}(T) \to \textrm{GL}_{2}(T_{O}(T))$. Once again, we have that $\textrm{Ker}(\pi) = Id$ and $\pi(\textrm{GL}_{2}(T)) \subset O_{2}(\mathbb{R})$.
Let $\sigma_{1}, \ldots \sigma_{N} \in T_{O}(T)$ be such that $\sigma_{i} = \textrm{hol}(S_{i})$, and let $\tau_{1}, \ldots \tau_{N} \in T_{O}(T)$ be such that $\tau_{i} = \textrm{hol}(T_{i})$.
The sets $\{\sigma_{1}, \ldots \sigma_{N}\}$ and $\{\tau_{1}, \ldots \tau_{N}\}$ are $\pi(\textrm{GL}_{2}(T))$-invariant, and both are the set of vertices of a regular $N$-gon with center at the origin. Both sets are invariant under a rotation of $\theta$ about the origin. However, due to the algebraic independence of $\ell_{1},\ell_{2},\ell_{3}$, no reflection that leaves $\{\sigma_{1}, \ldots \sigma_{N}\}$ invariant can also leave $\{\tau_{1}, \ldots \tau_{N}\}$ invariant. Thus, $C_{N} \cong \pi(\textrm{GL}_{2}(T)) \cong \textrm{GL}_{2}(T)$.

\paragraph{The case where N is odd.} Choose numbers $0 < \ell_{3} < \frac{1}{2}\ell_{2} <  \ell_{2} < \frac{1}{3} \ell_{1} < \ell_{1}$, such that $\ell_{1}$, $\ell_{2}$ and $\ell_{3}$ are algebraically independent. Construct a polygon as in the odd case of the proof of Proposition $3.1$, with the following difference: when constructing the square with sides of length $\ell_{2}$ that shares a side with the isosceles triangle, we require that the midpoints of the sides that the squares share \textit{do not} coincide. Let the distance between the midpoints be $\ell_{3}$. After identifying congruent parallel sides of the resulting polygon, we get a translation surface which we call $T$. In this case there are no new ideas to show that $\textrm{GL}_{2}(T) \cong C_{N}$. The proof is simply a restatement of the arguments found in the even case of Proposition $3.2$ and the odd case of Proposition $3.1$, and we omit it. $\Box$

\paragraph{Proof of Theorem 1.} The existence of translation surfaces whose Veech group is cyclic of order $2$ is well known (See the survey in  \cite{Mollsur} for a proof of this fact.) The remaining cases are treated in Propositions $3.1$ and $3.2$. $\Box$

\bibliography{finiteveech}
	\bibliographystyle{plain}

\end{document}